\documentclass[12pt]{article}

\usepackage{amsmath,amssymb,amsthm}
\usepackage[english]{babel}
\usepackage[latin1]{inputenc}
\usepackage{graphicx}
\usepackage[colorlinks,urlcolor=blue]{hyperref}

\newtheorem{proposition}{Proposition}
\newtheorem{theorem}{Theorem}
\newtheorem{lemma}{Lemma}
\newtheorem{remark}{Remark}

\newcommand{\EE}[1]{I\!\!E\left[#1\right]}
\newcommand{\NN}{I\!\!N}
\newcommand{\norm}[1]{\left\vert#1\right\vert}

\newcommand{\RR}{I\!\!R}

\newcommand{\prtl}{\partial}

\author{\small Haidar Al-Talibi\\
\small Department of Mathematik\\
\small Linnaeus University \\
\small Haidar.Al-Talibi@lnu.se\\
\and
\small Astrid Hilbert\\
\small Department of Mathematik\\
\small Linnaeus University\\
\small Astrid.Hilbert@lnu.se \\
\and
\small Vassily Kolokoltsov \\
\small Department of Statistics\\
\small University of Warwick\\
\small v.Kolokoltsov@warwick.ac.uk}
\title{Smoluchowski-Kramers Limit for a System \\
        Subject to a Mean-Field Drift}

\begin{document}

\parindent = 0pt

\maketitle

\begin{abstract}
We establish a scaling limit for autonomous stochastic Newton equations, the
solutions are often called nonlinear stochastic oscillators, where the
nonlinear drift includes a mean field term of McKean type and the driving
noise is Gaussian. Uniform convergence in $L^2$ sense is achieved by
applying $L^2$-type estimates and the Gronwall Theorem. The approximation is
also called Smoluchowski-Kramers limit and is a particular averaging
technique studied by Papanicolaou. It reveals an approximation of diffusions
with a mean-field contribution in the drift by stochastic nonlinear
oscillators with differentiable trajectories.
\end{abstract}
\section{Introduction}
In E. Nelson~\cite{ENelson} Brownian motion is constructed as a scaling
limit of a family of Ornstein-Uhlenbeck position processes which possess
differentiable sample paths by construction. The so-called ''Ornstein-Uhlenbeck theory of
Brownian motion'' constitutes a dynamical theory of Brownian motion. In
principle, this result goes back to~\cite{Smol} and~\cite{Kramers} and also
Chandrasekhar~\cite{Chand} studied this kind of limits. The second order
stochastic Newton equation is represented as a degenerate system of first
order It\^{o} equations. For more details see Nelson~\cite{ENelson} and the
references therein.

Degenerate diffusion processes which satisfy stochastic Newton equations
driven by Brownian motion have been investigated in many works
e.g.~\cite{AHK,AHZ, Narita1991}. For applications in the setting of diffusions
on manifold, see~\cite{Elworthy}. The broader class of Lévy processes also
includes processes with jumps and hence can explain more phenomena.
We proved a Smoluchowski-Kramers limit for stochastic Newton equations driven by
Lévy noise see~\cite{AlHiKo, Al} in the absence of mean-field terms. In the
PhD thesis of Zhang~\cite{Zhang} the reader finds references concerning
generalizations to an infinite dimensional setting. Also it should be mentioned
that the Smoluchowski-Kramers approximation carries over to nonlinear
stochastic oscillators driven by Fractional Brownian motion~\cite{AlHiKo2012, CKM}. We would
like to emphasize that the mathematical techniques and the mode of convergence in
the works mentioned above differ substantially.
\\
There are studies prior to this work considering the differential approximation
in the presence of mean field, one of these is~\cite{Narita1991}. In that paper,
a differentiable approximation for the stochastic Liénard equation with mean-field
of McKean type is achieved.

In the theory of differential games with a large number of players or in financial
mathematics in particular for pricing models for markets with a large number of
traders, stochastic differential equations with drift of mean-field type appear
naturally~\cite{UH}. Here the mean-field term contains no dissipation. In this paper we
pursue a scaling limit as treated in~\cite{AlHiKo, Al, ENelson} for stochastic
Newton equations driven by Brownian motion with a nonlinear drift term
$(\beta K)$, $\beta > 0$, where $K$ is globally Lipschitz continuous, moreover,
it depends only on the position process and a mean-field term in form of the
expectation value of the position process.
\section{Preliminaries}
In this section we recall a representation for It\^o processes with linear
dissipation and explain the scaling limit we are using.
\begin{proposition}
Consider the equation
\[
dY_t=-aY_t dt+f(t)dt+dB_{t},
\]
where $a>0$, $B_{t}$, $t\geq 0$,  is standard Brownian motion,  and
$f:\RR\to\RR$ is a measurable function such that
$\int_{0}^{\tau}e^{as}|f(s)|ds<\infty$.
Then we know that the solution exists and is given by
\begin{equation}\label{eq:sdelemma}
Y_t=e^{-at}x_0+\int_{0}^{t}e^{-a(t-s)}f(s)ds+\int_{0}^{t}e^{-a(t-s)}dB_{s},
\end{equation}
\end{proposition}
\noindent
for some initial value $x_0$. For a proof see for example~\cite{Protter2003}.
\\
A stochastic Newton equation is an It\^o equation which is second order in time,
hence is given by a system of a first order ordinary differential equation and a
stochastic differential equation. The dependent variables are traditionally
denoted by $x_t$ and $v_t$, respectively, whence  the independent variable
$t\geq 0$ is interpreted as time.
In physical models $x_t$ describes the position of a particle at time $t>0$.
It is assumed in this paper that the velocity $\frac{dx}{dt}=v$ exists and
satisfies the so called Langevin equation with a nonlinear drift depending on
the position and its expectation value.
\\
We introduce a scaling in form of a parameter $\gamma > 0$ for the stochastic
Newton equation with mean-field having the following form
\begin{eqnarray*}\label{eq:1}
dx_t&=& v_t\,dt\\
dv_t&=&-\gamma v_t dt + K(x_t-\EE{x_t})dt + \sqrt{\gamma} dB_t.
\end{eqnarray*}
Deviating from Narita, in our model the nonlinear part of the drift includes the
mean field term $\EE{x_t}$, moreover, it does not contain the marginal process
$v_t$ which is why this term does not explicitly contain any scaling parameter.
We emphasize that the deterministic system corresponding to this system
apart from the dissipation $\gamma v_t$ does not depend on the scaling
parameter. Moreover, for $ K(x)=-\frac{\prtl V}{\prtl x}$ we may define a
Hamiltonian $H(x,v)= \frac 12 v^2 + V(x)$, such that
$\frac{\prtl H}{\prtl v}= \frac{\prtl H}{\prtl v}= v$ and
$\frac{\prtl H}{\prtl v}=-\beta K(x)$. Averaging properties of such Hamiltonian
systems have been studied in \cite{AlK94, Li}.

For the solution of~($\ref{eq:1}$) we scale the time according
to $t'=\frac{1}{\gamma}t$, and we define the scaled process
$x^{\gamma}(t')=x(\gamma t')$ and $v^{ \gamma}(t')=\gamma v(\gamma t')$.
Moreover, we introduce a new Brownian motion
$\tilde{B}_{t'} \stackrel{d}{=} \frac{1}{\sqrt{\gamma}}B(\gamma t')$ for a 
new Brownian motion $\tilde{B}$.
Then $(x^{\gamma}(t'),v^{\gamma}(t'))$ satisfy the following stochastic
differential equation
\begin{eqnarray}
dx_{t'}^{\gamma}&=& v_{t'}^{\gamma}\,d{t'}\\\nonumber
dv_{t'}^{\gamma}&=&-\gamma^2 v_{t'}^{\gamma} dt'
                  +\gamma^2 K(x_{t'}^{\gamma}-\EE{x_{t'}^{\gamma}})dt'
                  +\gamma^2 d\tilde{B}_{t'},
\end{eqnarray}
where in a final step the Langevin equation has been multiplied by $\gamma$. 
For the sake of a short notation we replace $t'$ by $t$,
$\tilde{B}_t$, by $B_t$, and
set $\gamma^2=\beta$. Then the two differential equations combine
to the initial value problem:
\begin{eqnarray}\label{eq:dvtwithmeanfiled}
dx_t^{\beta}&=& v_t^{\beta}\,dt\\
dv_t^{\beta}&=&-\beta v_t^{\beta} dt+\beta K(x_t^{\beta}-\mu_t^\beta)dt+\beta dB_t,
\end{eqnarray}
with initial state $(x^{\beta}_0,v^{\beta}_0)=(x_0,v_0)$, where
$\mu_t^\beta=\EE{x_t^\beta}$ denotes the expectation, the possibly nonlinear
function $K$ satisfies a global Lipschitz condition, i.e.
$|K(x-y)|\leq \kappa|x-y|$ for $\kappa>0$, and $B_t$ is standard Brownian 
motion.
\section{Formulation of the main result}
The system we finally derived will be investigated in this work. When examining
the equation in the variable $v_t$, apparently, the drift has the time scale
$\beta t$ and $\beta B_t \stackrel{d}{=} \tilde B_{\beta^2 t}$ has the faster time scale $\beta^2 t$.
Thus, for $\beta$ tending to infinity the time is sent to infinity while the
Brownian paths performs arbitrarily fast oscillations around the trajectories of
the deterministic system given by the drift.
\\
Let us formulate our main result. In a first step we introduce the It\^o process
in $\RR$ solving the stochastic differential equation
\begin{equation}\label{eq:dytsmolu}
dy_t=K(y_t-\EE{y_t})+dB_t,
\end{equation}
with $y(0)=x_0$ and $K$ as in~($\ref{eq:dvtwithmeanfiled}$).
\begin{theorem}\label{thm1}
Let $(x_t^{\beta},v_t^{\beta})$ be the solution of~(\ref{eq:dvtwithmeanfiled}) 
and let $\Phi=(x_0,v_0)$ be any two-dimensional random vector independent of 
the Brownian motion $B_t$, $t\geq 0$, such that
\begin{equation}
\EE{|\Phi|^{2}}< M < \infty.
\end{equation}
Then
\[
\EE{\sup_{0\leq t\leq T}\left|x^{\beta}_t-y_t\right|^{2}}\to 0
\quad \mbox{ as }\beta\to\infty,
\]
for every $T<\infty$, where $y(t)$ is the solution of~(\ref{eq:dytsmolu}).
\end{theorem}
The idea of the proof is to use Gronwall lemma repeatedly to show that $\mu_t$ and $\EE[\sup_{0\leq t \leq T}x_t^\beta]$ are uniformly bounded on $[0,T]$. Before we proceed let us make the following important remarks.
\begin{remark}\label{remark1}
As stated in the introduction the function $K$ is supposed to be globally
Lipschitz continuous throughout the paper. Moreover, assume all conditions
stated in Theorem ~\ref{thm1} hold. Then by Arnold~\cite{LArnold}, and for
more recent results see~\cite{IkedaWatanabe}, the differential
equations~(\ref{eq:dvtwithmeanfiled}) and~(\ref{eq:dytsmolu}) have a
pathwise unique solution $\left(x^{\beta}(t),v^{\beta}(t)\right)$ with
initial states $\left(x_0,v_0\right)$ satisfying the moment estimates
\[
\EE{\sup_{0\leq t \leq T}\left|x_t^{\beta}\right|^{2}
            +\left|v_t^{\beta}\right|^{2}}\leq L_{\beta},
\]
for every $T<\infty$ with a constant $L_{\beta}>0$ depending on $\beta$ and
$T$.
\end{remark}
\begin{remark}
For the same assumptions as in Remark~\ref{remark1} we have for $y_t$
defined in (\ref{eq:dytsmolu}) with $y(0)=x_0$:
\[
\EE{\sup_{0\leq t \leq T} \left|y_t\right|^{2}}\leq L,
\]
for every $T<\infty$ and a constant $L>0$.
\end{remark}
\begin{remark}
The process $y_t$, $t\geq 0$, is defined in the same space as the coordinate
processes $x_t^\beta$, $t\geq 0$ in (\ref{eq:dvtwithmeanfiled}). In differential geometrical terms this means that
the driving Brownian motion $B_t$, $t \geq 0$, changed from the tangent
space $\RR$ where the coordinate process $v_t$, $t \geq 0$, is defined to
the manifold which trivially is $\RR$. Moreover, the nonlinear vector field
$K(x_t - \EE{x_t})$ changes from the cotangent space to the tangent space.
\end{remark}
Using Proposition~\ref{eq:sdelemma} the integral equation corresponding to the
first equation of~($\ref{eq:dvtwithmeanfiled}$) becomes
\begin{equation}\label{eq:xequation}
x_t^\beta = x_0 +v_0\int_{0}^{t}e^{-\beta s}ds
       +\beta\int_{0}^{t}\int_{0}^{s}e^{-\beta(s-u)}K(x_u^\beta-\mu_t^\beta)duds
       +\beta\int_{0}^{t}\int_{0}^{s}e^{-\beta(s-u)}dB_uds.
\end{equation}
We change the order in the double integral to obtain
\begin{eqnarray*}
x_t^{\beta}&=&x_0+I_0^\beta+\beta\int_{0}^{t}e^{\beta s}K\left(x_s^{\beta}-\mu_s^{\beta}\right)\int_{s}^{t}e^{-\beta u}duds+\beta\int_{0}^{t}e^{\beta s}\int_{s}^{t}e^{-\beta u}dudB_s,
\end{eqnarray*}
where $x_0=x^{\beta}(0)$ and
$I_0^\beta(t)=I^{\beta}(0,t)=\frac{v_0}{\beta}(1-e^{-\beta t})$.
Partial integration reveals
\begin{equation}\label{eq:equationx}
x_t^{\beta}=x_0+I_{0}^{\beta}(t)+I_1^{\beta}(t)+\int_{0}^{t}(1-e^{-\beta (t-s)})K(x_s^{\beta}-\mu_s^{\beta})ds+B_t,
\end{equation}
where $I_{1}^{\beta}(t)=-e^{-\beta t}\int_{0}^{t}e^{\beta s}dB_s.$ For the sake of
a short notation we sometimes drop the time parameter, i.e.
$I_1^\beta:= I_1^\beta(\cdot)$.
\section{Auxiliary estimates}
We evaluate the upper bound for the expectation value of each $I_{i}^{\beta}(t)$,
$i=0,1$, given above. For the deterministic integral $I_0^\beta$ the $k^{th}$
absolute moment, $k\in\NN $, is trivially estimated
\begin{eqnarray}\label{I1estimate}
 \EE{|I_0^\beta|^k}=\frac{1}{\beta^k}(1-e^{-\beta t})^k\leq\frac{1}{\beta^k}.
\end{eqnarray}
Moreover, for $I_1^\beta(t)$ we use It\^{o} isometry to have
\begin{equation}\label{I1estimatenew}
 \EE{|I_1^\beta(t)|^2}
 \le \EE{e^{-2\beta t}\left(\int_{0}^{t}e^{\beta t}dB_s\right)^2}
   =\frac{1}{2\beta} \left(1-e^{-2\beta t}\right)
  \le\frac{1}{2\beta}.
\end{equation}
Using the Lipschitz condition on $K$ we estimate $|x_t^\beta|$ in~($\ref{eq:equationx}$) by:
\begin{eqnarray}\label{eq:jensenineq}
 |x_t^\beta|\leq |x_0|+ \frac 1{\beta}+|I_1^\beta|+|B_t|
               +\kappa\int_{0}^{t} |x_s^\beta-\mu_s^\beta|ds.
\end{eqnarray}
By monotonicity of integration there holds $|\mu_t|\leq\EE{|x_t^\beta|}$. Taking
the expectation and inserting~(\ref{I1estimate}) and~(\ref{I1estimatenew})
reveals
\begin{eqnarray*}
 \EE{\norm{x_t^\beta}}
 &\leq \EE{| x_0 |}+\frac 1{\beta}+\EE{\norm{I_1^\beta}}+\EE{|B_t|}
                  +2\kappa\int_{0}^{t}\EE{| x_s^\beta|}ds\\
 &\leq  M^{\frac 12} +\frac{1}{\beta}+\frac{1}{\sqrt{2\beta}}+\sqrt{t}
                  +2\kappa\int_{0}^{t}\EE{| x_s^\beta |}ds,
\end{eqnarray*}
where we used the fact that $\EE{|x_0|^2}\le M $ and
$\EE{|B_t|}\le\sqrt{t}$. 
Next we consider the moment of the supremum of the square of the solution
$x^{\beta}(t)$,
for which we derive a bound uniform in the parameter $\beta$ and $t\in[0,T]$.
\begin{lemma}\label{lemma1}
For the same assumption as in Theorem~\ref{thm1}, let
$\left(x^{\beta}_t,v^{\beta}_t\right)$ be the solution
of~(\ref{eq:dvtwithmeanfiled}) with initial state $\Phi$. Then we have
\[
 \sup_{\beta>1}\EE{\sup_{0\leq t \leq T}\left|x^{\beta}(t)\right|^{2}}
 \leq H(T)
\]
for arbitrary fixed $T<\infty$, where $H(T)$ is a positive
constant independent of $\beta$ as $\beta$ tends to infinity.
\end{lemma}
\subsubsection*{Proof of Lemma~\ref{lemma1}}
Let $T<\infty$ be arbitrary but fixed. Consider any $t\in[0,T]$.
We return to equation~(\ref{eq:jensenineq}), square both sides, apply Jensen's
inequality before taking the supremum over the interval $[0,t]$, and estimate:
\begin{equation}\label{eq:estisup}
 \sup_{0\leq u\leq t}|x_t^\beta|^2
 \leq 5\left(|x_0|^2 \!+\! \frac 1{\beta^2} 
             \!+\! \sup_{0\leq u\leq t}|I_1^\beta|^2
         \!+\!\sup_{0\leq u\leq t}|B_u|^2
         \!+\! 2\kappa^2 T \int_{0}^{t}|x_s^\beta|^2
         \!+\!\EE{|x_s^\beta|^2} ds\right)
\end{equation}
where we have used that we have positive integrands and the estimate
$|\mu_t|^2\leq\EE{|x_t^\beta|^2}$. Since $B_t$ is a martingale,
Doob's inequality reveals $\EE{\sup_{0\leq u\leq t}|B_u|^2}\le 4\EE{|B_t|^2}$,
respectively, by inserting (7) and the recent estimates into
\begin{displaymath}
 \sup_{0\leq u\leq t}| I_1^\beta (u) |^2
 \le 2 \sup_{0\leq u\leq t}|B_u|^2
   + 2 \sup_{0\leq u\leq t}|B_u|^2 t\sup_{0\leq u\leq t} (1-e^{-2\beta u}).
\end{displaymath}

we obtain that for all $0\le t\le T$:
\begin{equation}\label{eq:estim}
\EE{\sup_{0\leq  u\leq t}|x_u^\beta|^2}
 \leq 5\left(\EE{|x_0|^2} + \frac{1}{\beta^2} + 9T + 8T^2
  +4\kappa^2 T\int_{0}^{t}\EE{\sup_{0\leq u\leq s}|x_u^\beta|^2}ds\right).
\end{equation}
Let us now turn to the core of the proof Lemma~\ref{lemma1}.
By the assumption on the initial values $x_0, v_0$ we have
\[
  \EE{|x_0|^{2}}< M < \infty.
\]
For $\beta>1$ and hence for $\beta\rightarrow\infty$ on a given interval $[0,T]$
we define the constant
\[
  H_0(T):= 5 M +5 + 45T + 40T^2
\]
which combines all additive constants in equation (\ref{eq:estim}) and rewrite
\begin{eqnarray*}
 \EE{\sup_{0\leq u \leq t}\left|x^{\beta}(t)\right|^{2}}
\leq \left[H_0(T)+20\kappa^2 T\int_{0}^{t}\EE{\sup_{0 \leq u \leq s}
             \left|x_u^{\beta}\right|^{2}}ds\right].
\end{eqnarray*}
Gronwall's lemma then reveals:
\begin{equation}\label{xestimate3}
 \EE{\sup_{0\leq u \leq t}\left|x^{\beta}(t)\right|^{2}}
    \leq  H(T),
\end{equation}
with $H(T)=H_0(T)e^{\theta T^2}$ and $\theta:=20\kappa^2$ where $\kappa$ is the
Lipschitz constant of $K$. The bound holds uniformly for $t\in[0,T]$ and for
$\beta$ sufficiently large, e.g. $\beta>1$, and hence for $\beta\rightarrow\infty$,
which completes the proof.
\\
Next we give the proof of the main result of this paper.
\section{Proof of Theorem~\ref{thm1}}
Under the same assumptions as in Theorem~\ref{thm1}. Let
$(x_t^{\beta}, v_{t}^{\beta})$ and $y_t$ be the solutions
of~($\ref{eq:dvtwithmeanfiled}$) and~($\ref{eq:dytsmolu}$) with initial states
$(x^{\beta}_0, v^{\beta}_0)=\Phi=(x_0,v_0)$ and $y_0=x_0$, respectively.
\\
Combining~($\ref{eq:dytsmolu}$) and~($\ref{eq:equationx}$) we have
\[
 x_t^{\beta}-y_t=\int_{0}^{t}K(x_s^{\beta}-\mu_s^{\beta})ds
          -\int_{0}^{t}K\left(y_s-\mu_s\right)ds+I_{0}^{\beta}(t)
               +I_{1}^{\beta}(t)+I_{2}^{\beta}(t),
\]
where $I_{i}^{\beta}(t), i=0,1$, are as before,
$I_{2}(t)=-e^{-\beta t}\int_{0}^{t}e^{\beta s}K(x_s^\beta-\mu_s^\beta)ds$,
and~$\mu_s=\EE{y_s}$. We exploit that $K$ is Lipschitz continuous with constant
$\kappa$ and rearrange the arguments in the norm, to get
\begin{eqnarray*}
 |x_t^{\beta}-y_t|
 &\leq |I_0^\beta|+|I_1^{\beta}|+|I_2^\beta|
                  +\int_{0}^{t}|K(x_s^\beta-\mu_s^\beta)-K(y_s-\mu_s)|ds\\
 &\leq\kappa\int_{0}^{t}|x_s^\beta-y_s|ds+\kappa\int_{0}^{t}|
 \mu_s^\beta-\mu_s|ds
     +\sum_{i=0}^{2}\left|I_{i}^{\beta}(t)\right|.
\end{eqnarray*}
By inserting the definition of $\mu_s^\beta$ and $\mu_s$ and by using Jensen's
inequality, in particular that
$|\mu_s^\beta - \mu_s|^2 \le \EE{\left|x_s^{\beta} - y_s\right|^2}$,
we estimate further
\begin{eqnarray*}
 |x_t^\beta-y_t|^2
 \leq 5\left(\kappa^2t\int_{0}^{t}|x_s^\beta-y_s|^2ds
      +\kappa^2t\int_{0}^{t}\EE{|x_s^\beta-y_s|^2}ds
      +\sum_{i=0}^{2}\left|I_{i}^{\beta}(t)\right|^2\right).
\end{eqnarray*}
Taking first the supremum and then the expectation value yields for
$0\leq t  \leq T<\infty$, that
\begin{displaymath}\label{eq:estimate5sup}
 \EE{\sup_{0 \leq u \leq t}\left|x_u^{\beta}-y_u\right|^{2}}
 \leq  D\int_{0}^{t} \EE{\sup_{0 \leq u \leq s}\left|x_u^{\beta}-y_u\right|^{2}}ds
      + 5\sum_{i=0}^{2}\EE{\sup_{0 \leq u \leq t} \left|I_{i}^{\beta}(u)\right|^2}
\end{displaymath}
with $D:= 10\kappa^2 T $. Direct calculation together with the previous bounds
(\ref{I1estimate}) and (\ref{I1estimatenew}) reveals that the last term on the
right hand side is uniformly bounded. To this end we proceed in the same way as
in (\ref{eq:estisup}) and (\ref{eq:estim}), in particular we use that $K$ is
Lipschitz continuous as well as Jensen's inequality to find
\begin{eqnarray*}
 \left|I_{2}^{\beta}(t)\right|^2
 &\le& t e^{-2\beta t}\int_{0}^{t}e^{2\beta s}(K(x_s^\beta-\mu_s^\beta))^2 ds
  \leq \kappa^2 t e^{-2\beta t}\int_{0}^{t}e^{2\beta s}
            |x_s^\beta-\mu_s^\beta|^2 ds\\
 &\le& 2\kappa^2 t e^{-2\beta t}\int_{0}^{t}e^{2\beta s}|x_s^\beta|^2 ds.
\end{eqnarray*}
Since the integrand is positive taking the supremum reveals the following estimate:
\begin{eqnarray*}
 \sup_{0\leq u \leq t}|I_2^\beta(u)|^2
 &\leq 2\kappa^2 t\sup_{0\leq u \leq t}|x_s^\beta|^2
                   \frac{1}{2\beta}(1-e^{-2\beta t})
   \leq\kappa^2 T \sup_{0\leq u \leq t}|x_u^\beta|^2 \frac{1}{\beta}.
\end{eqnarray*}
Taking expectation and inserting the bound (\ref{xestimate3}) we obtain:
\begin{equation}\label{eq:estimateI2supnew}
 \EE{\sup_{0\leq u \leq t}\left|I_{2}^{\beta}(u)\right|^{2}}
 \leq \frac{\kappa^2}{\beta} T H(T)e^{\theta T^2}.
\end{equation}
For $\beta$ sufficiently large, e.g. $\beta \ge 1$, we get  the following bound
\begin{eqnarray}\label{allIestimate}
\sum_{i=0}^{2}\EE{\sup_{0\leq u\leq t}\left|I_{i}^{\beta}(u)\right|^{2}}
\le  \frac{1}{\beta} \Lambda(T)
\end{eqnarray}
with $\Lambda(T) := 5 \left(\kappa^2 T H(T)e^{\theta T^2}+\frac{3}{2}\right)$
where we combined the bounds (\ref{I1estimate}), (\ref{I1estimatenew}), and
(\ref{eq:estimateI2supnew}).
By introducing the bound~(\ref{allIestimate}) into~($\ref{eq:estimate5sup}$) we arrive at the inequality
\begin{eqnarray*}
 \EE{\sup_{0 \leq u \leq t}\left|x_u^{\beta}-y_u\right|^{2}}
 \leq \frac 1\beta \Lambda(T) + D \int_{0}^{t}
   \EE{\sup_{0 \leq u \leq s}\left|x_u^{\beta}-y_u\right|^{2}}ds
\end{eqnarray*}
which allows to apply Gronwall's lemma a last time, to give
\begin{eqnarray*}
 \EE{\sup_{0 \leq u \leq t}\left|x_u^{\beta}-y_u\right|^{2}}
 &\leq \frac 1\beta \Lambda(T) e^{D T^2}
\end{eqnarray*}
with constant $D$ as defined in (\ref{eq:estimate5sup}) which holds uniformly on
the given time interval $[0,T]$ and does  not depend on $\beta$ for $\beta$
sufficiently large, while
$\frac{1}{\beta}$ tends to zero as $\beta$ tends to infinity.
This concludes the proof of Theorem~\ref{thm1}.
\section{Section}
\textbf{Acknowlwdgement}
The authors gratefully acknowledge Stig Larsson's kind comments on the first draft of
the paper and stimulating discussions with Francesco Russo and Sergio Albeverio.
The first author gratefully acknowledges the grant FOA12Magn-014, Kungl.
Vetenskapsakademien.

\end{document}